\documentclass{article}
\usepackage[%
journal=JST,    
lang=british,   
]{ems-journal}

\usepackage[margin=0pt]{subfig}

\numberwithin{equation}{section}

\newtheorem{theo}{Theorem}[section]
\newtheorem{rem}{Remark}[section]
\newtheorem{prop}{Proposition}[section]
\newtheorem{defi}{Definition}[section]
\newtheorem{lemm}{Lemma}[section]

\newcommand{\R}{{\mathbb{R}}}
\newcommand{\C}{{\mathbb{C}}}
\newcommand{\N}{{\mathbb{N}}}
\newcommand{\Z}{{\mathbb{Z}}}

\newcommand{\go}{\omega}
\newcommand{\gl}{\lambda}
\newcommand{\pa}{\partial}
\newcommand{\ha}{\frac{1}{2}}
\newcommand{\ra}{\rightarrow}
\newcommand{\gs}{\sigma}
\newcommand{\gd}{\delta}
\newcommand{\ga}{\alpha}
\newcommand{\gb}{\beta}

\newcommand{\BOPD}{pseudo-differential boundary operator}
\newcommand{\FIO}{Fourier Integral Operator}
\newcommand{\OPD}{pseudo-differential operator}

\begin{document}

\title{The spectrum of the Poincar\'e operator in an ellipsoid}
\titlemark{The spectrum of the Poincar\'e operator in an ellipsoid}



\emsauthor{1}{
	\givenname{Yves}
	\surname{Colin de Verdi\`ere}
	\mrid{50580}
	\orcid{0000-0001-7350-4662}}{Y.~Colin de Verdi\`ere}
\emsauthor{2}{
	\givenname{J\'er\'emie}
	\surname{Vidal}
	\orcid{0000-0002-3654-6633}}{J.~Vidal}

\Emsaffil{1}{
	\department{Institut Fourier}
	\organisation{Universit\'e Grenoble Alpes, CNRS},
	\zip{38000}
	\city{Grenoble}
	\country{France}
	\affemail{yves.colin-de-verdiere@univ-grenoble-alpes.fr}}
\Emsaffil{2}{
	\department{Laboratoire de G\'eologie de Lyon: Terre, Plan\`etes et Environnement}
	\organisation{CNRS, Universit\'e Lyon 1, ENS de Lyon},
	\city{Lyon}
	\country{France}
	\affemail{jeremie.vidal@ens-lyon.fr}
	}

\classification{35Q35}[53Z05, 76U60, 76B70]


\begin{abstract}
We study the spectrum of the Poincar\'e operator in triaxial ellipsoids subject to a constant rotation. 
As explained in the paper, this mathematical problem is interesting for many physical applications. 
It is known that the spectrum of this bounded self-adjoint operator is pure point with polynomial eigenvectors [Backus \& Rieutord, Phys. Rev. E \textbf{95} (2017), 053116]. 
We give two new proofs of this result. 
Moreover, we describe the large-degree asymptotics of the restriction of that operator to polynomial vector fields of fixed degrees. 
The main tool is the microlocal analysis of the partial differential equation satisfied by the orthogonal polynomials in ellipsoids.
This work also contains numerical calculations of these spectra, showing a very good agreement with the mathematical results. 
\end{abstract}

\maketitle

\section{Introduction}
Large-scale flows in natural objects (e.g. planetary liquid cores or stars) are often subject to global rotation.
A striking feature of such rotating flows is the ubiquitous presence of inertial waves (or modes in some  geometries).
These wave motions, which exist even without density effects for incompressible flows, are sustained by the Coriolis force \cite{greenspan1968theory}. 
If the rotating fluid has a non-zero and spatially uniform vorticity $ \vec{\Omega} \in \R^3$ (we define  $\go:= \lVert \vec{\Omega} \rVert$), these motions are in the simplest case small-amplitude perturbations governed by the linearised rotating Euler equation  and the incompressible condition
\begin{subequations}
\label{equ:tNSL} 
\begin{equation}
u_t  + \vec{\Omega}\wedge u= - \nabla \psi, \quad \mathrm{div} (u) = 0,
\tag{\theequation a--b}
\end{equation}
\end{subequations}
where the vector $u$ is the fluid velocity and the scalar $\psi$ is the pressure. 
Inertial modes can be excited by various mechanisms in natural objects, such as orbital (mechanical) forcings \cite{aldridge1969axisymmetric,rieutord2010viscous,lin2023resonant} or turbulent convection \cite{zhang1994coupling,lin2021triadic}. 
Moreover, inertial modes are often key in the dynamics of rapidly rotating fluids.
For instance, nonlinear couplings of inertial modes can sustain flow instabilities \cite{kerswell2002elliptical,vidal2023precession}, turbulence \cite{grannan2017tidally,le2019experimental} and, possibly, planetary (or stellar) magnetic fields through dynamo action \cite{reddy2018turbulent,vidal2018magnetic}.

\begin{figure}
    \centering
    \begin{tabular}{cc}
    \begin{tabular}{c}
    \includegraphics[width=0.4\textwidth]{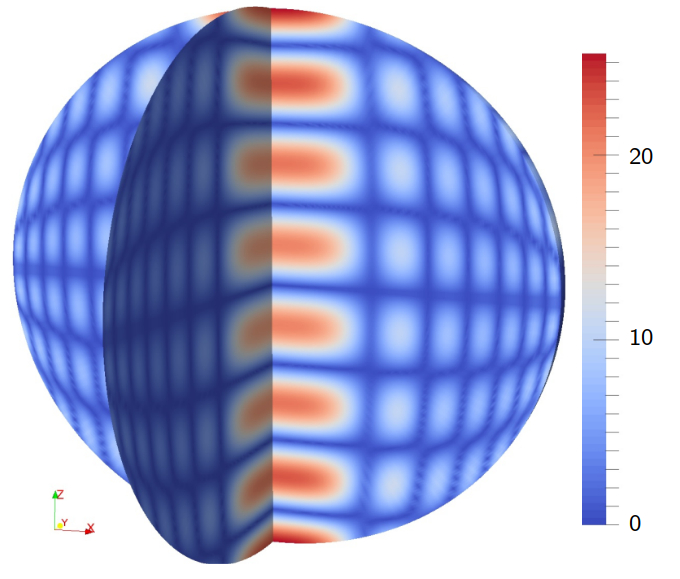}
    \end{tabular} &
    \begin{tabular}{c}
    \includegraphics[width=0.45\textwidth]{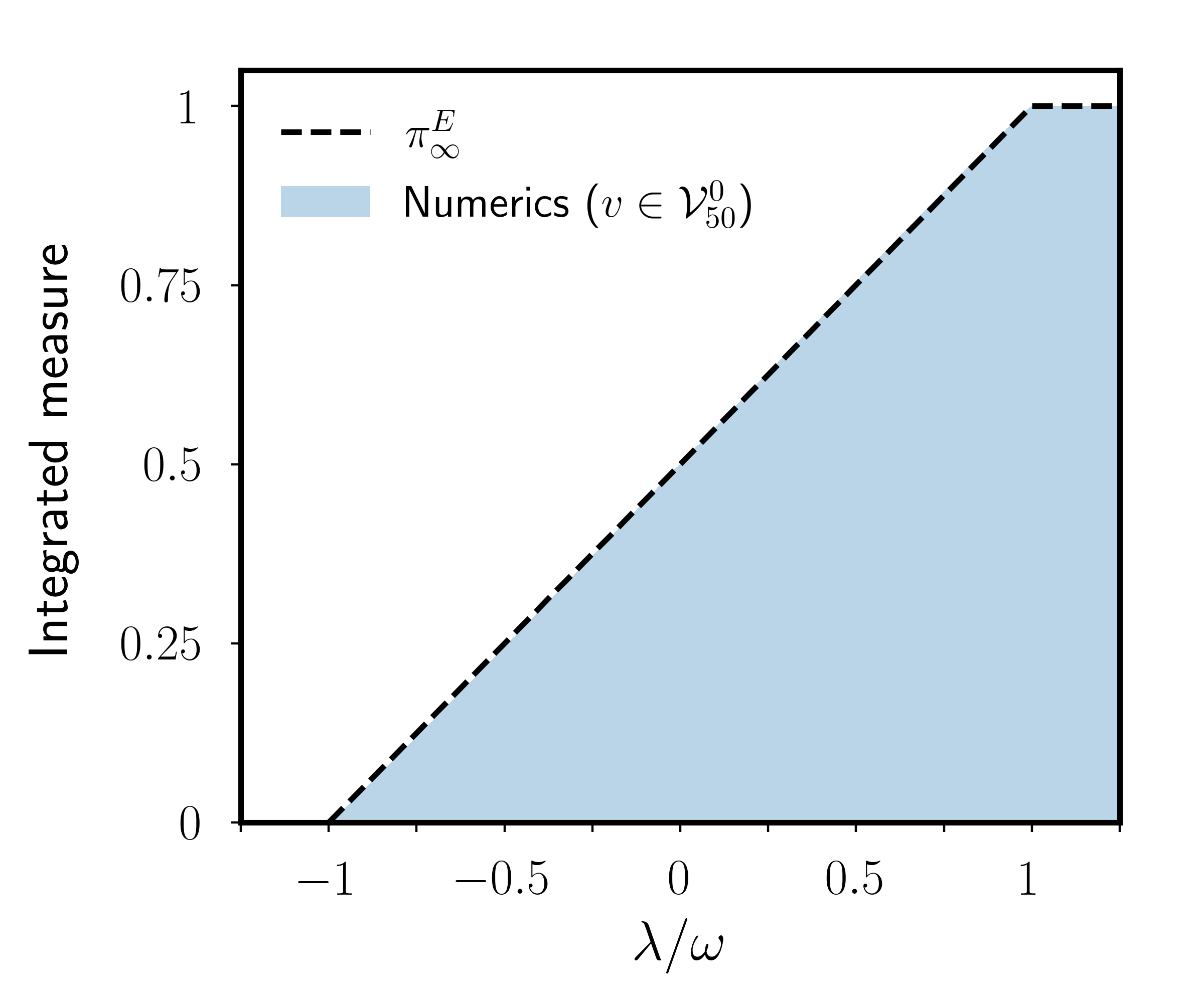}
    \end{tabular}\\
    (a) & (b) \\
    \end{tabular}
    \caption{(a) Velocity magnitude $\lVert v \rVert$ of an inertial mode with angular frequency $\lambda/\omega \simeq 0.5412$ in an axisymmetric ellipsoid ($A_1=A_2=1$, $A_3\simeq1.179$) with $\vec{\Omega} = (0,0,1)$. (b) Integrated probability measure of all the eigenvalues for $v \in {\mathcal V}_{50}^0$ in the round ball \cite{vidal2020compressible}, which  converges towards a uniform distribution.}
    \label{fig:intro}
\end{figure}

Owing to their considerable importance for planetary or astrophysical applications, inertial modes (and their related flows) are often studied in ellipsoidal geometries \cite{le2015flows}.
Indeed, rapidly rotating fluid masses are usually ellipsoidal at the leading order \cite{chandrasekhar1969ellipsoidal} (because of centrifugal forces and, possibly, tidal interactions due to orbital partners).
Let us give some definitions.
For $A_1,A_2, A_3 >~0$, we introduce the ellipsoid 
\begin{equation}
E=\{x=(x_1,x_2,x_3)\in \R^3 \mid A_1x_1^2+A_2x_2^2 +A_3x_3^2 \leq 1 \}.
\end{equation}
We equip $\R^3$ with the canonical Euclidean structure and the canonical orientation. 
We denote by $\mathcal{V}$ the Hilbert space of vector fields in $E$ whose coefficients are in $L^2(E, |dx|)$, where $|dx|$ is the Lebesgue measure, and by ${\mathcal{V}}^0$ the closed subspace orthogonal to the vector fields that are gradients of smooth functions. 
A smooth element in ${\mathcal{V}}^0$ is divergenceless and tangent to the boundary \cite[chapter 3]{Ga}.
Inertial modes are periodic solutions $u=\mathrm{e}^{i\gl t}v$ of equation \eqref{equ:tNSL}, where the complex-valued eigenvector $v\in {\mathcal{V}}^0$ is given by
\begin{subequations}
\label{equ:NSL} 
\begin{equation}
i\gl v + \vec{\Omega}\wedge v = -\nabla \psi, \quad {\mathrm{div}} (v) = 0,
\tag{\theequation a--b}
\end{equation}
\end{subequations}
together with the no-penetration condition saying that $v$ is tangent to  the boundary of the ellipsoid  $\partial E$. 
An example of a large-scale inertial mode in an axisymmetric ellipsoid is shown in figure \ref{fig:intro}~(a).
Even in this simple physical configuration, solving the inertial mode problem is very challenging from a mathematical viewpoint.
This is more clearly evidenced by considering the equation for the pressure (called the Poincar\'e equation after Cartan \cite{Ca}, who revisited Poincar\'e's paper \cite{Po})
\begin{subequations}
\label{eq:poincareeq}
\begin{equation}
    \lambda^2 \Delta  \psi = (\vec{\Omega} \cdot \nabla)^2 \psi, \quad \langle \nabla \psi \mid \vec{n} \rangle |_{\partial E} = \langle v \wedge \vec{\Omega} \mid \vec{n} \rangle |_{\partial E},
    \tag{\theequation a--b}
\end{equation}
\end{subequations}
where $\vec{n}$ is the unit vector normal to the boundary $\partial E$. 
The Poincar\'e equation is hyperbolic for $|\gl | < \omega$,  but, because of the  boundary condition, the inertial-mode problem is an ill-posed Cauchy problem \cite{rieutord2000wave}.
Explicit solutions in $\mathcal{V}^0$ have been found in spheroids with $A_1=A_2$ and $\vec{\Omega}=(0,0,1)$
since the pioneering work of Bryan \cite{bryan1889vi}, which admit (Cartesian) polynomial expressions \cite{kudlick1966transient,zhang2017theory}. 
Low-degree polynomial solutions in ${\mathcal{V}}^0$ have also been found in triaxial ellipsoids \cite{vantieghem2014inertial}, but there are no explicit solutions for the higher-degree modes. 
Actually, the inertial-mode spectrum is pure point in ellipsoids and the polynomial eigenvectors form a complete set \cite{B-R,ivers2017enumeration}.
The latter result could open new lines of research in fluid dynamics \cite{greenspan1968theory}, but the inertial-mode spectrum is still not well understood in an ellipsoid.

In this paper, we aim to better understand the properties of the inertial-mode spectrum in ellipsoids.
We first prove, using another mathematical route, that the inertial-mode spectrum is pure point with polynomial eigenvectors. 
The same result was initially proved in axisymmetric cases \cite{De}, and recently extended to triaxial geometries \cite{B-R,ivers2017enumeration}.
Then, we present new mathematical results on the asymptotic behaviour of the pure point spectrum in triaxial ellipsoids, which are successfully validated against numerical computations.
Note that the methods presented in this work can be extended to two-dimensional \cite{CdV-Li} or three-dimensional \cite{vidal2024igw,CdV-Vi-2} internal gravity waves.


\section{Poincar\'e operator in ellipsoids}
\subsection{Some generalities}
We denote by $\mathcal{P}_n$ the space of polynomial functions from $\R^3$ to $\C$ of degree less or equal to $n$, by
$\mathcal{V}_n$ the space of vector fields in $E$ whose coefficients are in $\mathcal{P}_n $ and by
$\mathcal{V}_n ^0=\mathcal{V}_n\cap \mathcal{V}^0 $
the subspace of  $ \mathcal{V}_n $ whose elements are the vector fields that are divergenceless and tangent to the boundary $\pa E$.

If we denote by $\Pi $ the orthogonal projection from $\mathcal{V} $ onto  $\mathcal{V}^0 $, called the Leray projector, we get the Poincar\'e operator $P $, which is a bounded self-adjoint operator on $\mathcal{V}^0$ defined by
\begin{equation}
    P v =i \Pi \Big ( \vec{\Omega}\wedge \Pi v \Big ).
\end{equation}
The spectrum of $P$ is the set of values of $\gl$ for which there are non-zero solutions of equation \eqref{equ:NSL}.
We have the following result in ellipsoids \cite{B-R}
\begin{theo}
\label{theo:A}
The spaces $\mathcal{V}_n ^0 $ are invariant by $P$. 
The spectrum of $P$ is the interval $[ -\go, \go   ]$ and is pure point. There is an orthonormal basis of $\mathcal{V}^0 $ consisting of eigenvectors of $P$ that have  polynomial coefficients.
\end{theo}
Note that this result is interesting for (at least) two reasons.
First, the fact that the spectrum is pure point shows that there are no attractors in the classical dynamics. 
Second, the fact that the eigenmodes are polynomial allows a good numerical calculus of the spectra.
 
Next, we introduce the spaces $\mathcal{W}_n =\Big ( \mathcal{V}_{n-1}^0 \Big )^\perp \cap \mathcal{V}_n^0 $ of dimension $d_n =n(n+2)$.
It follows from Theorem \ref{theo:A} and from the self-adjointness of $P$ that the spaces $\mathcal{W}_n $ are invariant by $P$. 
Our main result is about the asymptotic repartition of the eigenvalues $\mu_j^n$, for $j\in \{1,\dots,d_n\}$, of the restriction of $P$ to the spaces $\mathcal{W}_{n}$.
To this end, we define the probability measures $\pi_n $ by
 \begin{equation}
     \pi_n:=\frac{1}{d_n}\sum_{j=1}^{d_n} \delta (\mu_j^n ),
 \end{equation}
and we have the following theorem.
\begin{theo}
\label{theo:B}
As $n\ra \infty$, the measures $\pi_n$  converge weakly to a probability  measure $\pi_\infty $ of support $[ -\go, \go   ]$.
It means that, for every continuous function $f:\R\ra \R $, we have
\begin{equation}
\lim _{n\ra \infty }\frac{1}{d_n}\sum_{j=1}^{d_n} f(\mu_j^n ) =\int _\R f \, d\pi_\infty.
\end{equation}
\end{theo}
\begin{rem} A weaker result is obtained by looking not at a fixed degree $n$, but at the joint repartition of the eigenvalues $\mu_j^l$ for $j \in \{1,\dots,d_l\}$ and $l\leq n $. 
The limits are the same.
This is this last measure that is numerically computed below in \S\ref{sec:numericsasymptotic}. 
\end{rem}
Previous numerical computations suggest that the probability density is uniform in the round ball (see figure \ref{fig:intro}b),
but the properties of the spectral measure have not been  investigated in ellipsoids before our  work.
The recipe for the construction of $\pi_\infty$ will be given in Section \ref{sec:B}.
Note that the previous results could also be a good starting point for studying classical problems, such as other spectral asymptotics or control theory (e.g. see in section
\ref{sec:persp} below). 

\subsection{First proof of Theorem \ref{theo:A}}
We give here a first simple proof of Theorem \ref{theo:A}, which is quite close in spirit of previous proofs \cite{B-R,ivers2017enumeration} but without using dimension arguments. 
We first show that $\mathcal{V}_n=\mathcal{V}_n^0 \oplus \nabla \mathcal{P}_{n+1}$.
It follows from the Green formula that both spaces are orthogonal. Let us look at the orthogonal space of
$\nabla \mathcal{P}_ {n+1}$ in $\mathcal{V}_n$.
If $v=a\pa_1 +b\pa _2 +c\pa _3 $\footnote{We make use of differential geometry in considering a vector field as a differential operator. We also use the notations $\pa_j=\pa/ \pa x_j $.} is in $\mathcal{V}_n$ and orthogonal to $\nabla \mathcal{P}_ {n+1}$, we have for any polynomial $\phi \in \mathcal{P}_ {n+1}$
\begin{equation*}
\int _E \phi ~ \mathrm{div}(v) \, d{x_1} d{x_2} d{x_3} -\int _{\pa E} \phi  \langle v \mid \vec{n} \rangle \, d\sigma =0,
\end{equation*}
where
\begin{equation*}
\vec{n}= \Big(A_1^2{x_1}^2 + A_2^2 {x_2}^2 + A_3^2 {x_3}^2  \Big)^{-\ha} \Big( A_1{x_1} \pa_{1} + A_2{x_2} \pa_{2} + A_3{x_3} \pa_{3} \Big)
\end{equation*}
is the unit normal vector to $\pa E$.
By taking 
\begin{equation*}
    \phi = (A_1{x_1}^2+A_2{x_2}^2+A_3{x_3}^2-1) \, \mathrm{div}(v) \in \mathcal{P}_{n+1},
\end{equation*}
we get
\begin{equation*}
\int _E (A_1{x_1}^2+A_2{x_2}^2+A_3{x_3}^2-1) \, \mathrm{div }(v)^2 \, d{x_1} d{x_2} d{x_3}=0
 \end{equation*}
such that $\mathrm{div} (v) =0$.
Let us then take 
\begin{equation*}
\phi =aA_1{x_1}+bA_2{x_2} +cA_3{x_3} \in \mathcal{P}_{n+1},
\end{equation*}
we get
\begin{equation*}
\int_{\pa E} \Big({x_1}aA_1+{x_2}bA_2 +{x_3}cA_3 \Big)^2 \Big(A_1^2 {x_1}^2+A_2^2 {x_2}^2+A_3^2 {x_3}^2 \Big)^{-\ha} \, d\sigma =0
\end{equation*}
such that $v$ is tangent to $\pa E$.
The spaces $\mathcal{V}_n $ are invariant by $\vec{\Omega}~\wedge $. 
It follows that the spaces $\mathcal{V}_n^0$ are invariant
by the Poincar\'e operator $P=i\Pi \Big ( \vec{\Omega}\wedge \Pi \Big)$. 

We have to show that $\oplus_{n\in\N } {\mathcal{V}_n}^0$ is dense in $\mathcal{V}^0 $ for the $L^2$ topology.
By the Stone-Weierstra\ss~ Theorem, $\oplus_{n\in\N } \mathcal{V}_n $  is dense in $\mathcal{V}$. 
Moreover, the space 
\begin{equation*}
K:= ( \oplus_{n\in\N } \mathcal{V}_n^0 ) \bigoplus  ( \oplus_{n\in\N } \nabla \mathcal{P}_{n+1} )
\end{equation*}
is dense in $\mathcal{V} = \mathcal{V}_0 \oplus \nabla H^1$. 
The respective densities of each term in $K $ then follows.
We have an orthonormal basis of $L^2$ eigenmodes, and the spectrum is pure point. 
Finally, it remains to prove that the spectrum is the full interval. 
This follows from the methods of Theorem 2.1 in  \cite{CdV-b}.
$P $ is inside $E$ a pseudo-differential operator, and  the computation of the principal symbol $p$  of $P$ (see Section \ref{sec:B}) shows that the  image of the eigenvalues of $p$ is the  interval $[-\go, \go]$. 
On the other hand, we have $\lVert P \rVert \leq \go $. 
Hence, the full spectrum is $[-\go, \go ]$. 
This holds for any bounded smooth domain. 
For ellipsoids, it also results from Theorems \ref{theo:B} and \ref{theo:dens}.

\subsection{Further spectral properties}
\begin{prop}The numbers $\pm \go $ are not eigenvalues of $P$ in ellipsoids.
\end{prop}
\begin{proof}
The above result was proved in \cite{B-R,ivers2017enumeration}, but we give below an alternative proof for $-\go $. 
We assume that $v\in \mathcal{V}_0 $ satisfies $\Pi (\Omega \wedge v) =i\go v $. 
Then, by virtue of the Pythagorean theorem, we have $\vec{\Omega} \wedge v \in \mathcal{V}_0$ and $\vec{\Omega} \wedge v =i\go v $.
It follows that $v_3=0$ and
$v_1=iv_2$. Then, using $\mathrm{div}(v)=0$,  we see that $v_1 $ is harmonic in $(x_1,x_2)$. 
Moreover, $(n_1 -in_2)v_1=0$ on $\pa E$  implies $v_1=0$ on $\pa E$. 
To finish, observe that on each level set $z=c$, $v_1$ is harmonic and vanishes on the boundary.
So, we have $v_1=iv_2=0 $ everywhere. 
Actually, the previous argument is valid for any smooth bounded domain in $\R^3$.
\end{proof}

The subset of geostrophic modes, which are invariant along the rotation axis, often plays
an important role in the dynamics of rapidly rotating flows \cite{greenspan1968theory}.
A vector field $v \in \mathcal{V}_0 $ is \textit{geostrophic} if $v$ is in the kernel of the Poincar\'e operator. 
Without loss of generality, we can assume that $\vec{\Omega }=(0,0,1)$ and
\begin{equation*}
E = \{ x\in \R^3 \mid Ax_1^2+Bx_2^2+x_3^2 + 2C x_1x_3 +2Dx_2x_3 \leq 1 \}.
\end{equation*}
This equation for $E$ is more convenient in the study of geostrophic fields. 
We have used the scaling invariance and a rotation around the $x_3-$axis to make the coefficient of $x_3^2$ equal to $1$, and to have a vanishing coefficient for $x_1x_2$. 
We have the following result in ellipsoids (initially proved in \cite{ivers2017enumeration}).
\begin{prop}
There exists exactly one geostrophic field  in $\mathcal{W}_n$ for $n$ odd and no geostrophic field in $\mathcal{W}_n$ for $n$ even.
\end{prop}
\begin{proof}
Let us denote by $F$ the projection of $E$ onto the $(x_1,x_2)-$plane. For each $(x_1, x_2) $ in the interior of $F$, there exist
two points $(x_1,x_2, x_3^\pm (x_1,x_2))$ in $\pa E $.
From $\vec{\Omega} \wedge v= -\nabla \psi $, we get that  $v=(v_1,v_2,v_3)$ satisfies
\begin{equation*}
   v_1 = -\pa_1 \psi, \quad v_2 = \pa_2 \psi, \quad \pa_3 \psi=0,
\end{equation*}
and the boundary condition
\begin{equation*}
\Big( Ax_1+ Cx_3^\pm \Big) \, v_1 + \Big( Bx_2+ Dx_3^\pm \Big) \, v_2 + \Big( x_3^\pm + Cx_1 + Dx_2 \Big) \, v_3 = 0.
\end{equation*}
From $\pa_3 \psi=0$, we get that $\psi$ is independent of $x_3$.
From the expressions of $v_1$ and $v_2$ in terms of $\psi$, we get that they are independent of $x_3$ too.
Using $\mathrm{div}(v) = 0$, we get that $v_3$ is also independent of $x_3$.
Te four functions $[v_1,v_2,v_3,p]$ are thus independent of $x_3$.
This is called the Taylor-Proudman theorem \cite{greenspan1968theory}. 
Eliminating $v_3$ from the previous two equations and replacing $v_1$ and $v_2$ in terms of $p$, we get that $p$ satisfies a differential equation $Vp=0$ where the coefficients of $V$ are linear in $(x_1,x_2)$ such that
\begin{equation*}
V = \Big( CDx_1- (B-D^2)x_2 \Big) \, \pa_1 + \Big( (A-C^2)x_1 -CDx_2 \Big) \, \pa_2.
\end{equation*}
If we write $V=M(x) \, \pa_x$, we see that $\mathrm{Trace}(M)=0$. Moreover, the determinant of $M$ is given by
\begin{equation*}
\delta = AB-BC^2-AD^2,
\end{equation*}
which is also the determinant of the quadratic form defining $E$.
Thus, we have $\delta >0$. 
We can then find a basis of $\R^2$, where $V=a(u\pa_v -v\pa_u)$  with $a\ne 0$. 
This implies easily that  there is a unique (up to dilation) non-zero quadratic form $Q$ with $VQ=0$ and no linear forms $L$ with $VL=0$: $Q=u^2+v^2$. 
All polynomials $f$ are then clearly of the form $f=F(Q)$ where $F$ is a polynomial.
\end{proof}

\section{Orthogonal polynomials, Weyl law and a conceptual proof of Theorem \ref{theo:A}}
\subsection{Orthogonal polynomials in Euclidean balls}
Let us denote by $\mathcal{E}_n $ the space of polynomials of degree $n$ that are orthogonal to all polynomials of degree $\leq n-1$ in $L^2 (B, |dx|)$, where $B$ is the Euclidean ball of unit radius in $\R^3$.
The following result, due to Appell and Kamp\'e de F\'eriet \cite{A-K}, is proved in \cite[section 5.2]{D-X}:
\begin{theo}
\label{theo:ortho}
The spaces $\mathcal{E}_n$ are the eigenspaces of the operator $L$, which is called here the Legendre operator, defined by
\begin{equation*}
L= -\sum_{i=1}^3 \pa^2_i + \sum _{i,j=1}^3\pa_i x_i x_j \pa_j + \frac{9}{4}
\end{equation*}
with the eigenvalues $(n+3/2)^2$ of multiplicity $d_n^1={(n+1)(n+2)}/{2}$.
\end{theo}
\begin{rem} We can define an operator $L$ in any dimension by a similar formula. 
In dimension $1$, the  Legendre polynomials
are eigenfunctions of the operator $L=-\pa_x^2+ \pa_x x^2 \pa _x $. 
They are orthogonal polynomials on $L^2([-1,1],|dx|)$. 
\end{rem}

\begin{proof}
We give a simple proof of Theorem \ref{theo:ortho}.
For each $n\in\N$,  $L$ acts on  the space $\mathcal{P}_n $ of polynomials of degree $\leq n$. 
This action of $L$ is triangular.
If we decompose $\mathcal{P}_n $ into a direct sum of homogeneous polynomials 
\begin{equation*}
    \mathcal{P}_n = \bigoplus_{k=0}^n \, \mathcal{H}_k,
\end{equation*}
where $\mathcal{H}_k$ is the space of polynomials homogeneous of degree $k$, then we have 
\begin{equation*}
L \Big(\sum_{l=0}^k  h_l \Big)= (k+3/2)^2 h_k + r_k
\end{equation*}
with $\deg(r_k )\leq k-1$.
Hence, the eigenvalues of $L$ restricted to $\mathcal{P}_n $ are the numbers $(k+3/2)^2$ where $k\in\{0,\cdots,n\}$,
with the eigenspaces $\mathcal{H}_k \oplus \mathcal{R}_k$ with $\mathcal{R}_k \subset \mathcal{P}_{k-1}$. 
Let us show that $L$ is symmetric on $C^2(\bar{B})$, and hence on each $\mathcal{P}_n $. 
There is a cancellation of boundary terms coming from both parts of $L$. 
Let us rewrite \[ L =-\Delta -\mathcal{L}^\star \mathcal{L}, \]
with $\mathcal{L}= r\frac{d}{dr}$ and where $\mathcal{L}^\star $ is the formal adjoint of $\mathcal{L}$. 
By virtue of the Green-Riemann formula, we have
\begin{equation*}
\int_B (u\Delta v - v\Delta u) \, |dx| =\int_{\pa B} \Big( u\frac{\pa v}{\pa n} -v\frac{\pa u}{\pa n} \Big) \, |d\gs|.
\end{equation*}
On the other hand, we have also
\begin{equation*}
\int_B \Big(u \mathcal{L}^\star \mathcal{L} v  - v \mathcal{L}^\star \mathcal{L} u \Big) \, |dx| =-\int _{\pa B } \Big( u\frac{\pa v}{\pa r}-v\frac{\pa u}{\pa r} \Big) \, |d\gs|.
\end{equation*}
Both boundary terms cancel out since ${\pa}/{\pa r } = {\pa }/{\pa n}$ on $\pa B$.

It follows that the eigenspaces of $L$ are exactly given by the orthogonal polynomials.
The operator $L$ with domain $\oplus_{n\in \N} \, \mathcal{E}_n $ is essentially self-adjoint. It will be useful in particular in Section \ref{sec:micro}
to keep the notation $L$ for the differential operator, and to denote by $\widehat{L}$ its closure. 
Note that, if $u$ belongs to the Sobolev space $H^2(B)$, then $u$ belongs to the domain $ D(\widehat{L})$ of $\widehat{L}$.
\end{proof}

\subsection{Weyl law}
The principal symbol of $L$, denoted by $\Lambda $,  is given by
\begin{equation*}
\Lambda (x,\xi )=\lVert \xi \rVert^2 - \langle x \mid \xi \rangle^2. 
\end{equation*}
We see that $L$ is elliptic in the interior of $B$, but not on the boundary $\pa B$. 
The characteristic manifold is the co-normal bundle to the boundary.

The pull-back of $\Lambda $ onto the Euclidean sphere $S^3\subset \R^3_x \oplus \R_z $ by the orthogonal projection $(x,z)\ra x$ 
is the dual metric of the standard metric on $S^3$ (see Appendix \ref{ss:bic}). 
Therefore, the operator $L$ is very similar to the Laplacian on $S^3$ and the eigenfunctions similar to  the  spherical harmonics. 
However, the pull-back to $S^3$ of the Lebesgue measure on $B$ vanishes on the equator, and is not the canonical measure on $S^3$.
In fact, by looking at orthonormal polynomials in the ball with respect to the measure $(1-r^2)^{-\ha} |dx| $, we could make a similar analysis leading to the usual spherical harmonics (precisely the spherical harmonics that are even under the change $z\ra -z$). 

Let us look at the Weyl formula:
\begin{theo}\label{theo:weyl}
The eigenvalues counting function $N(\gl) $ of $L$ satisfies
\begin{equation*}
N(\gl) \underset{\gl\ra \infty }{\sim} \frac{\gl^{3/2}}{6},
\end{equation*}
where the notation $\sim $ means that the ratio goes to $1$ as $\gl \ra \infty $.
This expression coincides with  the phase-space volume calculated with respect to the Liouville measure such as
\begin{equation*}
\frac{1}{(2\pi)^3} \mathrm{Vol}(\{(x,\xi)\in T^\star B \mid \Lambda (x,\xi) \leq \gl\}) = \frac{\lambda^{3/2}}{6}.
\end{equation*}
\end{theo}
The Weyl formula can be easily derived from the explicit expression of the eigenvalues, which is
\begin{equation*}
N(\gl)= \sum _{n+3/2\leq \sqrt{\gl} } \frac{(n+1)(n+2)}{2}.
\end{equation*}
On the other hand, the calculus of the phase-space volume is a simple exercise. 
Theorem \ref{theo:weyl}  will be useful later in order to control the boundary effects. 
In the case of ellipsoids, we get similar  results by replacing $L$ by  $L_E$ given by
\begin{equation*}
L_E=-\sum _{i=1}^3 \frac{1}{A_i} \pa_i^2 +\sum_{i,j=1}^3 \pa _i x_i x_j \pa_j +\frac{9}{4}.
\end{equation*}
This is proved by using the affine diffeomorphism  $\Phi: B \ra E $ defined by
\begin{equation*}
\Phi (x_1,x_2,x_3)= \Big( \frac{x_1}{\sqrt{A_1}}, \frac{x_2}{\sqrt{A_2}},\frac{x_3}{\sqrt{A_3}} \Big),
\end{equation*}
and remarking that $\Phi $ transforms (i) the Lebesgue measure into a multiple of the Lebesgue measure and (ii) polynomials of degree $n$ into polynomials of degree $n$.
Note that $\Phi $ also transforms divergenceless vector fields that are tangent to the boundary in $B$ to vector fields with the same properties in $E$.

\subsection{A conceptual proof of Theorem \ref{theo:A}}
The Leray projector is the orthogonal projector on vector fields $L^2$-orthogonal to the space of gradient of smooth functions.
The operator $\sum_{i=1}^3 (1/A_i) \, \pa_i^2$ and the dilation operator $\mathcal{L}$ preserve the latter space.
$\mathcal{L}$ is called the \textit{dilation operator}, because this is the infinitesimal generator of the group of homotheties. 
The adjoint of $\mathcal{L}$ is $\mathcal{L}^\star = -\mathcal{L}-3 $. 
The Legendre operator given by
\begin{equation*}
L_E = -\sum _{i=1}^3 \frac{1}{A_i} \, \pa_i^2 - \mathcal{L}^\star \mathcal{L} + \frac{9}{4}
\end{equation*}
is self-adjoint, and it preserves the space of gradients. 
It implies that the Legendre operator $L_E$ commutes with
$\Pi$. 
This holds formally for any domain, but ${ L_E}$ is a well-defined symmetric operator only on ellipsoids. 
Note then that any operator with constant coefficients, such as $\vec{\Omega} \, \wedge$, commutes with $L_E$. Hence, the Poincar\'e operator  commutes with $L_E$.
This gives a proof of Theorem \ref{theo:A}  using only the spectral theory of $ L_E$. 

\section{The microlocal Weyl law}
\label{sec:micro}
This is the most technical part of this paper.
We use the \BOPD~calculus of Boutet de Monvel on manifolds with boundary \cite{Bo}, see also \cite{Gr} and Appendix \ref{app:lbm}, as well as the construction of a parametrix for the "wave equation" $u_{tt}+Lu=0 $ using \FIO s coming from \cite{Ho,D-G}.

Let $d_n^1:= \dim (\mathcal{E}_n)= (n+1)(n+2)/2 $. 
We have the following result:
\begin{theo}\label{theo:micro}
Let $A$ be a self-adjoint \BOPD~ of degree $0$ in $B$ of principal symbol $\gs(A)$.
Let us denote by $(\phi_j^n)$ with $j \in \{ 1,\cdots , d_n^1 \}$ an orthonormal basis of $\mathcal{E}_n$. 
In the limit of large $n$, we have 
\begin{equation*}
\lim_{n \ra +\infty} \frac{1}{d_n^1}\sum_{j=1}^{d_n^1} \langle A \phi_j^n \mid {\phi_j^n}\rangle = \frac{6}{(2\pi)^{3} }\int _{\Lambda \leq 1 }\gs (A) \, |dx d\xi|.
\end{equation*}
\end{theo}
\begin{proof}
We use the method of the papers \cite{CdV-a,We} that are  inspired from H\"ormander \cite{Ho}, see also \cite[section 2]{D-G}.
We have two main difficulties, namely (i) we have to work with a manifold with boundary and (ii) $L$ is not elliptic at the boundary.
We will make first an assumption on $A$ avoiding both difficulties, and then make approximations of $A$.
\end{proof}

\subsection{The case where $A$ is a "nice" \BOPD}
\label{ss:nice} 

If $A$ is a \BOPD, $\mathrm{WF}'(A)$ is the conic support of the full symbol of the \OPD~part of $A$. 
The conical set $Z$ is the set of points in the phase space $T^\star B \setminus 0$ such that the Hamiltonian flow of $\Lambda $ hits $\pa B $ at a characteristic point, namely
a point in  $N^\star \pa B $ (see Appendix \ref{ss:bic}).
We have
\begin{equation*}
Z=\{ (x,\xi )\in T^\star B \setminus 0 \mid \xi = sx ~\mathrm{for~some~} s \in \R \}.
\end{equation*}
In the following, it will be important to make a difference between $L$ as a differential operator acting in $\R^3$, and the self-adjoint operator on $L^2(B)$ denoted by $\widehat{L}$. 

\begin{defi}
\label{def:nice} 
A \BOPD ~$A$  is "nice" if it satisfies
\begin{equation*}
\mathrm{WF}'(A) \cap \Big( Z \cup \Big( \pa B \times \R^3 \setminus 0 \Big) \Big) = \emptyset .
\end{equation*}
\end{defi}

Recall that, if $K:\R \ra L(\mathcal{H})$, we can sometimes define the \textit{distributional trace} as the Schwartz distribution defined by
\begin{equation*}
\langle \mathrm{Trace}(K) \mid \phi \rangle = \mathrm{ Trace}\int_\R K(t)\phi(t) \, dt,
\end{equation*}
where $\phi$ is a test function, and the trace in the right-hand side is the usual trace of trace class operator.
The singularities of the distributional trace
\begin{equation*}
Z: t \ra Z(t):= \mathrm{Trace} \Big( \mathrm{e}^{-it\sqrt{\widehat{L}}}A \Big)=\sum_{n=0}^\infty a_n \mathrm{e}^{-it(n+3/2)}
\end{equation*}
with $a_n =\sum _{j=1}^{d_n^1} \langle A\phi _j^n \mid \phi_j^n  \rangle$, determine the asymptotics of the sequence $(a_n)_{n\in \N}$ (see Lemma \ref{lemm:tauber}).
 
We define $\sqrt{\widehat{L}} $ using the functional calculus for positive self-adjoint operators.
However, $L$ is not elliptic and $\sqrt{\widehat{L}} $ is not a  \BOPD~at characteristic points. 
It is more convenient to start from $\cos \Big( t\sqrt{\widehat{L}} \Big)$.
The solution of the wave equation
\begin{equation*}
\square u:= u_{tt} + \widehat{L}u =0, \quad u(0)=u_0, \quad u_t(0)=0
\end{equation*}
is $u(t) = u_0 \cos \Big( t\sqrt{\widehat{L}} \Big)$. 
If
\begin{equation*}
C: t\ra C(t):= \mathrm{Trace} \Big(\cos \Big(t\sqrt{\widehat{L}} \Big) A \Big),
\end{equation*}
then we have $Z = 2 H C$ where $H$ is the "Hilbert" projector defined by 
\begin{equation*}
\mathcal{F}(Hf)(\tau) = \mathrm{Heaviside}(-\tau) \mathcal{F}(f)(\tau)
\end{equation*}
where $\mathcal{F}(f)(\tau)=\int_\R \mathrm{e}^{-it\tau}f(t) \, dt$ is the Fourier transform of a function $f$. 
Hence, the singularities of $Z$ can be deduced from those of $C $. 

Let us denote by $\mathcal{E}:=\{ (x,\xi) \in T^\star \R^3 \mid \Lambda(x,\xi )>0 \} $ (the "elliptic" domain), by $\phi_t$ with $t\in \R$ the Hamiltonian flow of $\sqrt{\Lambda}$, and by
\begin{equation*}
\mathcal{L}_\pm := \{(t,\tau; x,\xi;y,\eta ) \mid (x,\xi)=\phi_{\pm t} (y,\eta ) \mid \tau \pm \sqrt{\Lambda}(x,\xi) =0,~(y,\eta )\in \mathcal{E} \}
\end{equation*}
the Lagrangian submanifolds of $T^\star (\R_t \times \R^3_x \times \R^3_y )$ associated with the flow $\phi_t$. 

\begin{lemm} If $A$ is a nice \OPD~ of degree $0$ (see the definition  \ref{def:nice}), the Schwartz kernel of $\cos \Big( t\sqrt{\widehat{L}} \Big) A$ is the restriction to $\R\times B \times B$ of a \FIO~ of degree $0$ associated with the union of the Lagrangian manifolds $\mathcal{L}_\pm$ modulo an operator $R(t)$ whose trace is a smooth function of $t$. 
\end{lemm}
\begin{proof}
In $\R^3$, the equation
\begin{equation*}
\square u =0, \quad u(0)=Au_0, \quad \dot{u}(0)=0,
\end{equation*}
admits a \FIO~parametrix $U(t)$ as in \cite[section 1]{D-G}. 
This is possible because $\mathcal{L}_\pm \subset T^\star \R \times \mathcal{E}\times \mathcal{E}$.
It means that $\square U(t)$ has a smooth kernel on $\R \times \R^3 \times \R^3$, and that $U(0)= \mathrm{A}$ and $U_t(0)=0$ modulo operators with smooth kernels.
We define then, for $u_0\in L^2 (B)$, $V(t)u_0= (U(t) u_0)_{|{{B}}}$.
The operator $V(t)$ is a parametrix for $\cos (t\sqrt{\widehat{L}}) A$:
it means that, if  $u_0 \in D(\widehat{L})$, $V(t)u_0\in D(\widehat{L})$, $S:=\square V(t)$ has a smooth kernel on $\R\times B \times B$ and that
$U(0)= \mathrm{A}$ and $ U_t(0)=0  $ modulo  operators with smooth kernels. 
If $u_0 \in D(\widehat{L})$, then $Au_0$ belongs to the Sobolev space  $ H^2(B )$. It follows from the continuity properties of FIO's that, for all $t\in \R$, $V(t)u_0 $ is also in $H^2 (B)$ and hence in the domain of $\widehat{L}$. 
The other properties follow from the properties of $U(t)$. 
If $\cos \Big( t\sqrt{\widehat{L}} \Big)A= V(t) +R(t)$, we have 
$\square R(t)= -S (t)$ with $R(0)$ and $\pa_t{R}(0)$ that are smoothing.
We get
\begin{equation*}
R(t) = \mathrm{e}^{it \sqrt{\widehat{L}}}R_1 + \mathrm{e}^{-it \sqrt{\widehat{L}}}R_2 +\int_0^t \mathrm{e}^{i(t-s)\sqrt{\widehat{L}}}R_3(s) \, ds +
\int_0^t \mathrm{e}^{-i(t-s)\sqrt{\widehat{L}}} R_4(s) \,ds,
\end{equation*}
where the different $R_j$ have smooth Schwartz kernels. 
It follows that they belong to the domains of all powers of $\widehat{L}$.
Hence, using the fact that $\mathrm{e}^{it \sqrt{\widehat{L}}}\widehat{L}^{-N} $ is a trace class for all $N \geq 2 $, we get the smoothness of $t\ra \mathrm{Trace}(R(t))$.
\end{proof}

The bicharacteristic flow lifts to the geodesic flow on $S^3$, which is simply periodic of period $2\pi$: the orbits starting outside $Z$ are projections of great circles of $S^3$, hence ellipses tangent to $\pa B $ (see Appendix \ref{ss:bic}).
This implies that there are no periodic orbits of period smaller than $2\pi $ in $\mathcal{L}_\pm $. 
Using the argument of \cite[corollary 1.2]{D-G}, we deduce from the calculi of wave-front sets that the distribution $C$ is only singular at the points $2\pi n$ with $n\in \Z$.
The antiperiodic distribution 
$Z= 2H C $ is also singular only at the points  $2\pi n$ with $n\in  \Z$.

The analysis provided in \cite{Ho,D-G} shows that the Fourier transform of $\rho Z $, with $\rho \in C_o^\infty (]-2\pi, 2\pi [)$ and $\rho(0)=1$, is a symbol $\tau \ra \gs (\tau )$ of degree $2$ and principal part
\begin{equation*}
\gs_2(\tau ) = \frac{6 \tau^2 }{(2\pi)^3} \int _{\Lambda\leq 1} \gs(A) \, |dx d\xi|.
\end{equation*}
Theorem \ref{theo:micro} in the case where $A$ is nice follows then from Lemma \ref{lemm:tauber}. 

\subsection{The general case}
\label{ss:generalcase}
Let us give some $\varepsilon >0$ and rewrite $A$ as $A=A_\varepsilon + R_\varepsilon$ in the following way.
We choose a \BOPD ~$Q_\varepsilon = \chi_\varepsilon +\eta_\varepsilon$ of degree $0$ so that, denoting by $d$ the Euclidean distances, 
\begin{enumerate}
\item the map $\chi_\varepsilon: \bar{B}\ra [0,1]$ is smooth, equal to $1$ if $d(x,\pa B) \leq \varepsilon / 2$ and vanishes if $d(x,\pa B)\geq \varepsilon$;
\item $\eta_\varepsilon$ is a compactly supported \OPD ~ in $B$ whose full symbol is $1$ in the cone $\{ (x,\xi) \mid d((x,\xi ), Z ) \leq \varepsilon/2 \lVert \xi \rVert ~\mathrm{and}~ d(x,\pa B) \geq \varepsilon /2 \} $ and $0$ in the cone $\{ (x,\xi) \mid d((x,\xi),Z) \geq \varepsilon \lVert \xi \rVert ~\mathrm{and}~ d(x,\pa B) \leq \varepsilon / 4 \}$.
\end{enumerate}

Let us decompose $A$ by putting $R_\varepsilon = Q_\varepsilon^\star A Q_\varepsilon$, and see that  $A - R_\varepsilon $ is nice. 
We have $A = P + G$, where $G $ is a Green operator (see Appendix \ref{app:lbm}). 
We have then
\begin{equation*}
A - R_\varepsilon= ( P-Q_\varepsilon ^\star P Q_\varepsilon ) + ( G - Q_\varepsilon^\star G  Q_\varepsilon ).
\end{equation*}
The first term $P - Q_\varepsilon ^\star P Q_\varepsilon$ is nice. 
Then, we need to prove that $B:=G - Q_\varepsilon^\star G Q_\varepsilon$ is smoothing.
We can write
\begin{equation*}
B = (\mathrm{Id} - Q_\varepsilon^\star )G +Q_\varepsilon^\star G (\mathrm{Id}-Q_\varepsilon),
\end{equation*}
and both terms are smoothing as follows from Proposition 10.11 in \cite{Gr} (see the precise statement in Appendix \ref{app:lbm}). 
We have to estimate $r_n(\varepsilon) :=\sum_{j=1}^{d_n^1} \langle R_\varepsilon \phi_j^n \mid \phi_j^n \rangle $. 
For any $\phi \in L^2$, we have
\begin{equation*}
| \langle R_\varepsilon \phi \mid \phi \rangle | \leq C \langle Q_\varepsilon^\star Q_\varepsilon  \phi \mid \phi \rangle
\end{equation*}
with $C= \lVert A \rVert $. 
Moreover, we have
\begin{equation*}
\sum_{j=1}^{d_n^1}\langle Q_\varepsilon^\star Q_\varepsilon \phi_j^n \mid \phi_j^n \rangle = d_n^1 - \sum_{j=1}^{d_n^1} \langle (\mathrm{Id} - Q_\varepsilon^\star Q_\varepsilon) \phi_j^n  \mid \phi_j^n \rangle.
\end{equation*}
The last term can be evaluated because $ \mathrm{Id}- Q_\varepsilon^\star Q_\varepsilon$ is nice, and we get
\begin{equation*}
\lim_{n\ra \infty } \frac{1}{d_n^1} \sum_{j=1}^{d_n^1} \langle \Big( \mathrm{Id} - Q_\varepsilon^\star Q_\varepsilon \Big)  \phi_j^n \mid \phi_j^n \rangle =
\frac{6}{(2\pi)^3} \int_{\Lambda \leq 1} (1 - \gs (Q_\varepsilon)^2) \, |dx d\xi|.
\end{equation*}
The integral $\int_{\Lambda \leq 1} \gs (Q_\varepsilon)^2  |d x d\xi|$ tends to 0 as $\varepsilon \ra 0$. 
It follows that 
\begin{equation*}
    \lim\limits_{n\ra \infty} \Big( \frac{1}{d_n^1} \Big) r_n (\varepsilon) \ra 0 \quad \mathrm{as} \, \varepsilon \ra 0.
\end{equation*}

\section{A scalar version of Theorem \ref{theo:B}}
We will prove the following theorem:
\begin{theo}
\label{theo:scalar}
Let $A$ be a self-adjoint \BOPD~of degree $0$ in the Euclidean ball $B$, which commutes with the operator $\widehat{L}$ and of principal symbol $\gs (A)$.
Let us denote by $\mu_j^n$ with $j\in \{1,\cdots,d_n^1\}$ the eigenvalues of $A$ restricted to $\mathcal{E}_n$. 
Then, when $n \ra \infty $,  the probability measures
\begin{equation}
\pi_n :=\frac{1}{d_n^1} \sum_{j=1}^{d_n^1} \gd \Big( \mu_j^n \Big)
\end{equation}
converge weakly to a probability measure $\pi_\infty $ that is defined as follows.
For any continuous function $f:\R \ra \R$, we have
\begin{equation}
\int_\R f \, d\pi_\infty 
= \frac{6}{(2\pi)^3} \int _{\Lambda\leq 1}f(\gs(A)) \, |dx d\xi|.
\end{equation}
\end{theo}
We can apply Theorem \ref{theo:micro} to the operators $A^N$ with $N\in\N$.
This gives 
\begin{equation*}
\langle A^N \phi_j^n \mid \phi_j^n \rangle = \Big( \mu_j^n \Big) ^N. 
\end{equation*}
Then, we get the asymptotics for $f$ a polynomial function.
The answer for a continuous function $f$  is given by uniform approximation of $f$ by polynomial functions. 

\section{The proof of Theorem \ref{theo:B}}
\label{sec:B}
\begin{lemm}\label{lemm:piopd}
The operator $\Pi $ is a \BOPD~in $B$ belonging to the Boutet de Monvel calculus (see Appendix \ref{app:lbm}).
The operator-valued symbol of $\Pi $ at the point $(x,\xi) \in T^\star B\setminus 0$ is the orthogonal projection on the hyperplane $\ker (\xi )$.
\end{lemm}
\begin{proof}
We use the isomorphism between vector fields and $2$-forms given by $V\ra \iota (V) dx_1 \wedge dx_2 \wedge dx_3$ where $\iota$ is the inner product. 
The image of the divergenceless vector fields tangent to the boundary is the closed $2$-forms $\alpha$ satisfying $j^\star (\alpha )=0$ where $j:\pa B \ra \R^3$ is the embedding.
The relative cohomology $H^2(B,\pa B )$ vanishes:
by Poincar\'e duality  $H^2(B,\pa B )$ is isomorphic to  $H_1(B)$ that vanishes. 
This implies that 
the Hodge Laplacian on $2$-forms with the relative boundary conditions is invertible and the inverse is a \BOPD ~in $B$.
The projector $\Pi$ is then given by $\Pi = d d^\star \Delta ^{-1}$, which is clearly also a \BOPD.
The symbol can be easily computed.
The symbol of the divergence is the linear form $\xi$, and the symbol of $\Pi$ is the orthogonal projection onto the hyperplane $\ker(\xi)$. 
\end{proof}

The spectrum of $P$ is the union of the eigenvalues of $P$ on $\mathcal{V}_0$, which are the eigenvalues $\gl_j ^n$ with $ n\in \N$ and $1\leq j \leq d_n$, and the eigenvalue $0$ with eigenspace $\mathcal{V}_0^\perp$. 
We can  compute the principal symbol $\gs({P})$ of the Poincar\'e operator
 using Lemma \ref{lemm:piopd}
and the composition rules of symbols.  If we  denote by $\phi \in [0,\pi ]$ the angle between $\vec{\Omega}$ and $\xi $,  the eigenvalues of
$\gs({P})(x,\xi)$ are $(0,\Lambda_1:=\go  \cos \phi ,-\Lambda_1 )$. If $\vec{\Omega}=(0,0,\go)$, 
this can be written as
\begin{equation*}
\Big( 0, \go \frac{\xi_3}{ \lVert \xi \rVert},- \go \frac{\xi_3}{\lVert \xi \rVert} \Big).
\end{equation*}
We see that these eigenvalues are simple if $\xi_3 \ne 0 $ where all eigenvalues vanish.

We will now prove Theorem \ref{theo:B}.
Let us recall that the Poincar\'e operator is given by $P=i \Pi (\vec{\Omega} \wedge \Pi)$
and that $P$ commutes with the extension to vector fields of the Legendre operator. We denoted by
$\mathcal{E}_n = \ker (L-n(n+3) \mathrm{Id})$ and are interested in the eigenvalues of $P$  restricted to $\mathcal{F}_n = \mathcal{E}_n \otimes \C^3$, which is a complex-valued vector space.  
We split these eigenvalues into the eigenvalues $0$ corresponding to $\ker \Pi _{|\mathcal{F}_n}$, which is of dimension $d_{n+1}^1$, and the other eigenvalues $\mu_j^n$ with $j \in \{1,\cdots,d_n\}$ and $d_n =3d_n^1 -d_{n+1}^1 = n(n+2)$.
We define
\begin{equation*}
\pi_n :=\frac{1}{d_n}\sum _{j=1}^{d_n} \gd \Big( \mu_j^n \Big),
\end{equation*}
and we are interested in the large $n$ limit of $\pi_n$.
We use Theorem \ref{theo:scalar} and replace $\gs(A)$ by the matrix-valued symbol $\gs(P)$ and $f(\gs (A))$ by
$\mathrm{Trace}f(\gs (P))$, such that we get
\begin{equation*}
\lim_{n\ra \infty}\int_\R f \, d\tilde{\pi}_n = \frac{6}{(2\pi)^3} \int_{\Lambda \leq 1} \mathrm{Trace}f(P) \, |dx d\xi|
\end{equation*}
where we have defined
\begin{equation*}
\tilde{\pi}_n:= \frac{1}{d_n^1} \left( d_{n+1}^1 \gd(0) + \sum  _{j=1}^{d_n} \gd \Big( \mu_j^n \Big) \right).
\end{equation*}
Then, we obtain
\begin{equation*}
f(0) + \frac{1}{d_n^1} \lim_{n\ra \infty}\sum_{j=1}^n f(\mu_j^n)  = \frac{6f(0)}{(2\pi)^3} \int_{\Lambda \leq 1 } \, |dxd\xi| +
\frac{6}{(2\pi)^3} \int_{\Lambda \leq 1} \Big( f(\Lambda_1) + f(-\Lambda_1) \Big) \, |dx d\xi|.
\end{equation*}
Recalling that $ \int_{\Lambda \leq 1} |dxd\xi | = (2\pi)^3/6 $ and $d_n\sim 2d_n^1$, we get
\begin{equation*}
\lim_{n\ra \infty} \int_\R f \, d\pi_n = \frac{3}{(2\pi)^3}\int _{\Lambda \leq 1 } \Big( f(\Lambda_1) + f(-\Lambda_1) \Big) \, |dxd\xi|.
\end{equation*}
Moreover, since $\Lambda_1(-\xi)=-\Lambda_1(\xi )$ and $\Lambda $ is even in $\xi$, we get
\begin{equation*}
\lim_{n\ra \infty }\int_\R f \, d\pi_n = \frac{6}{(2\pi)^3} \int_{\Lambda \leq 1 } f(\Lambda_1) \, |dxd\xi|.
\end{equation*}
Adapting to the ellipsoids, we get that the measure  $\pi_\infty^{E} $ is defined by
\begin{equation*}
\int_\R f \, d\pi_\infty^E = \frac{6}{(2\pi )^3}\int_{\gs(L_E)\leq 1} f(\Lambda_1) \, |dx d\xi|.
\end{equation*}
In other words, we have
\begin{equation}
\label{equ:piinfini}
\int_{-\infty}^t \, d\pi_\infty^E = \frac{6}{(2\pi )^3} \mathrm{Vol} \Big( \{(x,\xi)\in E \times \R^3 | \gs(L_E)(x,\xi)\leq 1 ~\mathrm{and}~ \go \xi_3\leq t \lVert \xi \rVert \} \Big).
\end{equation}

\section{Some properties of $\pi_\infty^E$}
\label{sec:numericsasymptotic}
\subsection{Asymptotic formula}
\begin{theo} \label{theo:dens} 
The spectra $(\mu_j^n)_{j=1}^{d_n}$ are symmetric with respect to $0$.
The measures $\pi_\infty^E $ admit densities $f_E:= d\pi_\infty^E/ du$ that are even, non-vanishing and  analytic in $[-\go,\go ]$. 
They have a jump at $\pm \go$.
Assuming that $\go=1$ for simplicity, we have
\begin{equation}
\label{equ:spectreasymp}
\forall u\in\R,~  \int_{-\infty}^u d\pi_\infty^E = \frac{1}{4\pi} \mathrm{Area}(C_u^E\cap S^2 )
\end{equation}
where $C_u^E$ is the cone of $\R^3$ defined by $C_u^E:=\{ \xi \mid \Lambda_E (\xi)\leq u \}$ where  
\begin{equation}
\label{equ:lambda}
\Lambda_E (\xi) = \go \cos \left( \vec{\Omega }, A_E (\xi ) \right) 
\end{equation}
and with $A_E(\xi) = \Big( \sqrt{A_1}\xi_1,\sqrt{A_2}\xi_2, \sqrt{A_3}\xi_3 \Big) $. 
\end{theo}
\begin{proof}
We will make the calculation using the coordinate system in $E$ given by the linear diffeomorphism 
$B_E : B \ra E $ defined by 
\begin{equation*}
B_E (x)= \Big (\frac{x_1}{\sqrt{A_1}}, \frac{x_2}{\sqrt{A_2}}, \frac{x_3}{\sqrt{A_3}} \Big),
\end{equation*}
and the induced canonical transformation
$\Phi^\star :T^\star B \ra T^\star E$ defined by 
$\Phi^\star (x,\xi)=(B_E(x),A_E(\xi))$. 
The pullback of the operator ${L}_E$ is ${L}$ while the pullback of $\omega \xi_3 / \lVert \xi \rVert$ is given by equation \eqref{equ:lambda}. 
The result is then obtained by calculating integral \eqref{equ:piinfini} as follows.
We first evaluate the integral in $x $ with $\xi $ fixed and, then, we calculate the $\xi$-integral in polar coordinates.
With more details, it works like that:
we need to compute
$I:= \int_{\Lambda (x,\xi)\leq 1} f(\xi ) \, |dx d \xi |$ with $f$ homogeneous of degree $0$. 
First, we compute 
$\mathrm{Vol} ( \{x\in B \mid \Lambda (x,\xi)\leq 1 \} )= F(\rVert \xi \rVert)$.
Using polar coordinates, we obtain
\begin{equation*}
I = \int _{\R^3} F(\lVert \xi \rVert)f(\xi ) \, |d\xi| = \int _0^\infty F(r)r^2 \, dr \times \int _{S^2}f(\gs) \, |d\gs| = C \int _{S^2}f(\gs ) \, |d\gs|,
\end{equation*}
where the constant $C$ is calculated by taking $f=1$. 
\end{proof}


\begin{figure}
\subfloat[]{\includegraphics[width=0.46\textwidth]{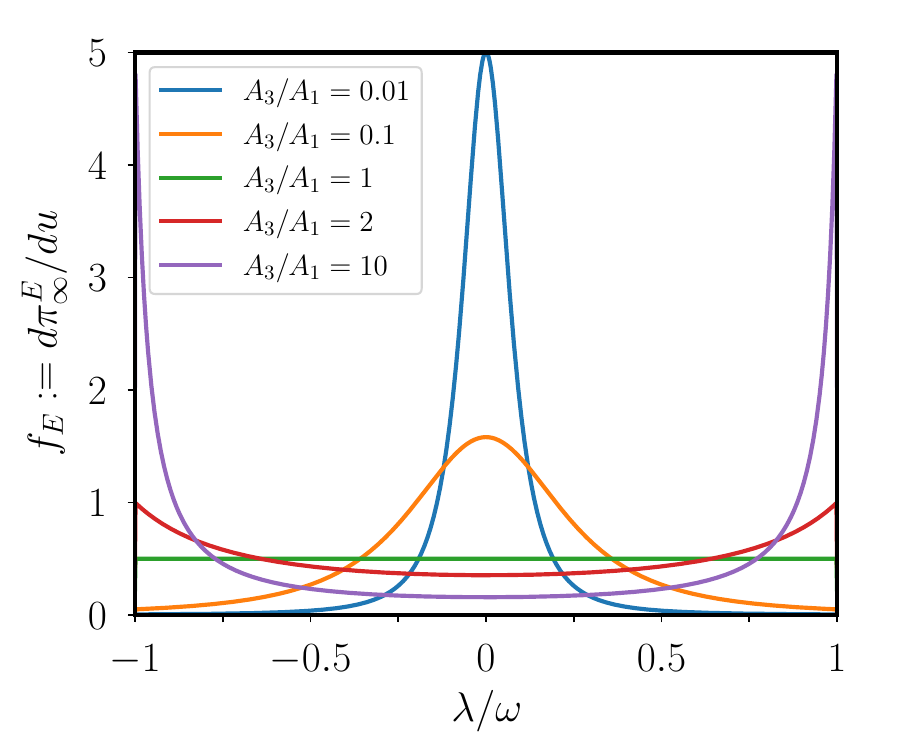}}
\qquad
\subfloat[]{\includegraphics[width=0.46\textwidth]{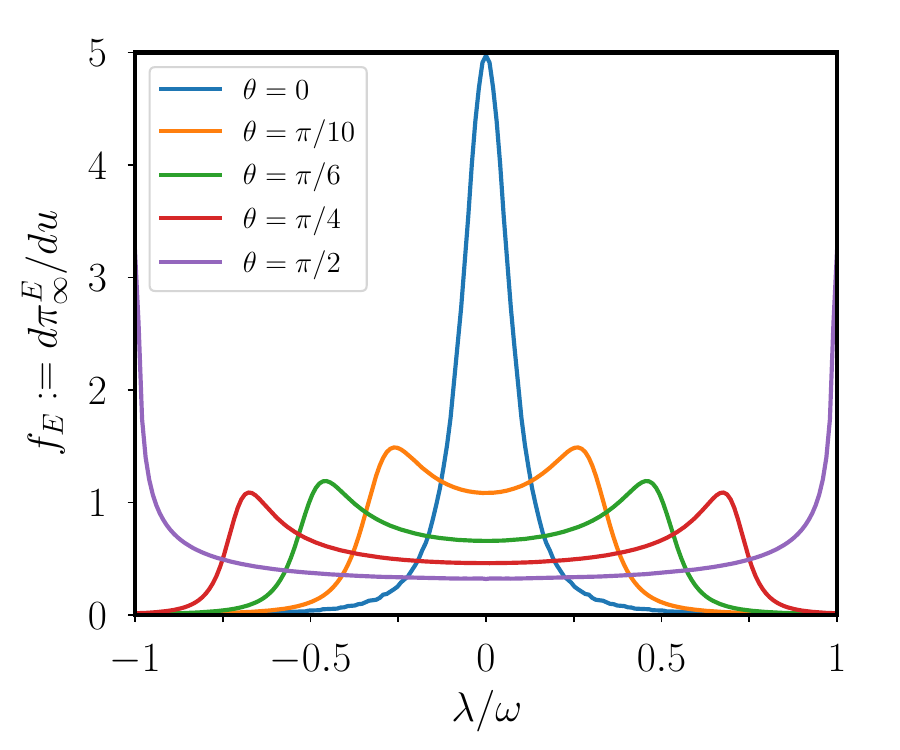}}
\caption{Non-vanishing density $f_E$ for eigenvalues $-1 < u < 1$ with $u := \lambda/\omega$ and $\omega=1$, obtained from Theorem \ref{theo:dens}. \textbf{(a)} Formula \eqref{eq:fEintegrable} for integrable cases in axisymmetric ellipsoids with $A_1=A_2=1$ and $\vec{\Omega }=(0,0,1)$. \textbf{(b)} Numerical evaluation of formula \eqref{equ:spectreasymp} for a tilted vector $\vec{\Omega }=(\sin \theta, 0, \cos \theta)$ in an axisymmetric ellipsoid $A_1=A_2=1$ and $A_3=0.01$.}
\label{fig:integrable}
\end{figure}

If $E$ is an ellipsoid of revolution with $A_1=A_2$ and $\vec{\Omega }=(0,0,1)$, the density of the measure $\pi_\infty^E$, defined as $f_E := {d\pi_\infty^E}/{du}$, is given by
\begin{subequations}
\label{eq:fEintegrable}
\begin{equation}
f_E = \chi_{[-1, 1]} (u) \frac{a}{2 [ u^2+ a(1-u^2) ]^{3/2} }, \quad \chi_{[-1, 1 ]} (u) = \begin{cases}
1 ~\text{if}~ |u| \leq 1, \\
0 ~\text{if}~ |u| > 1,
\end{cases}
\tag{\theequation a--b}
\end{equation}
\end{subequations}
where we have introduced $a={A_3/A_1}$. 
Formula \eqref{eq:fEintegrable} is illustrated in figure \ref{fig:integrable}~(a).
The density is uniform in $[-1,1]$ for the round ball, but we observe that the density becomes strongly non-uniform in flattened or elongated axisymmetric ellipsoids.
The density is indeed peaked near zero when $A_3/A_1 \to 0$, whereas the density becomes maximum near the edge of the spectrum $|\lambda| \to \omega$ when $A_3/A_1 \to \infty$.
Low-frequency inertial modes (known as quasi-geostrophic modes, since they are almost invariant along the rotation axis \cite{maffei2017characterization}) are thus favoured in prolate ellipsoids.
On the contrary, inertial modes have preferentially high frequencies in flattened ellipsoids.
For other configurations, formula \eqref{equ:spectreasymp} is not integrable (it generally involves elliptic integrals). 
Yet, the area intersection can be evaluated numerically (see figure \ref{fig:integrable}~b for tilted rotation axes).

\subsection{Numerical validation}
\begin{figure}[t]
\subfloat[Triaxial ellipsoid $A_1=1, \, A_2\simeq 1.235, \, A_3\simeq 2.04$ with $\theta=0$]{
\centering
        \begin{tabular}{cc}
        \includegraphics[width=0.46\textwidth]{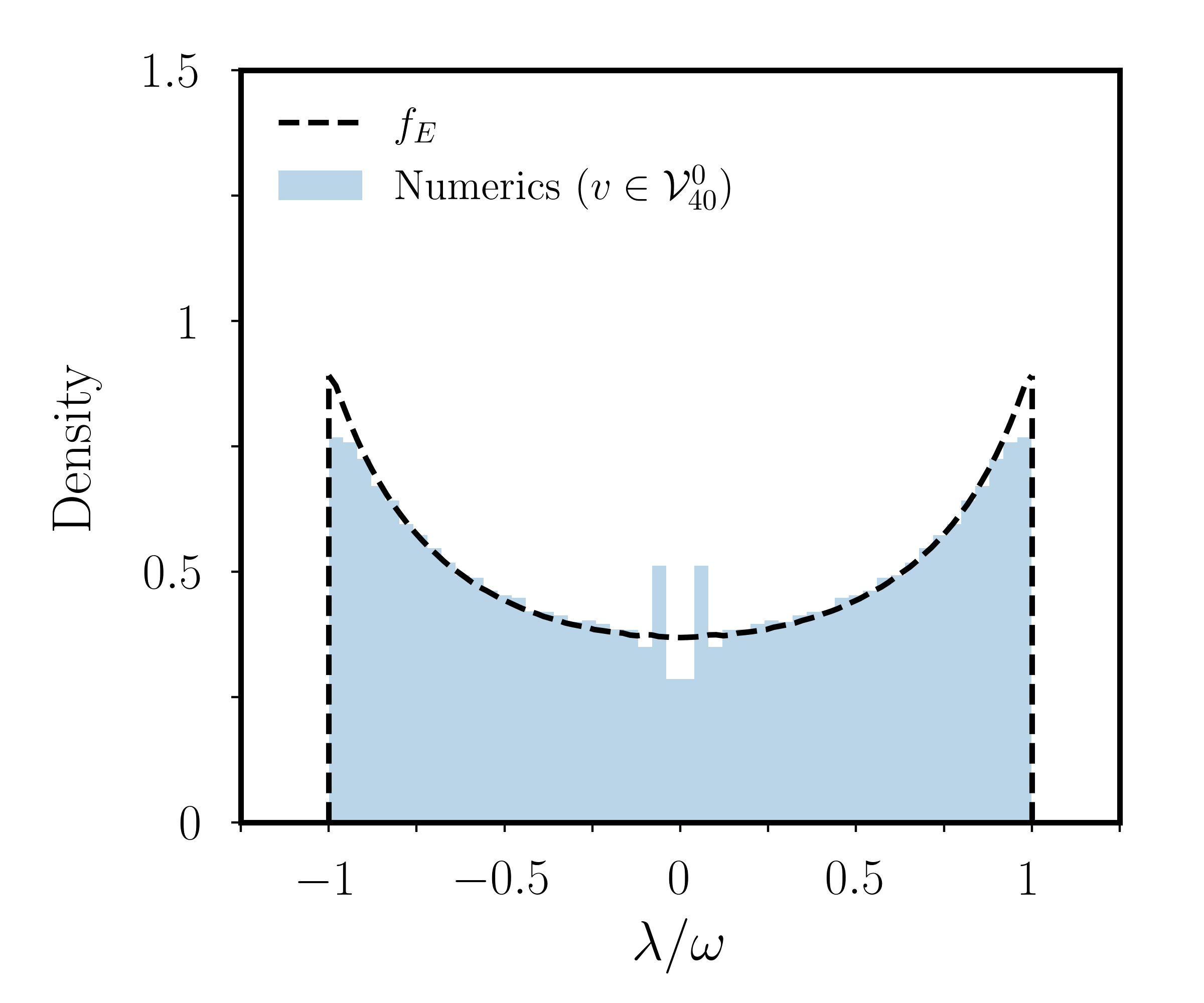}  &
        \includegraphics[width=0.46\textwidth]{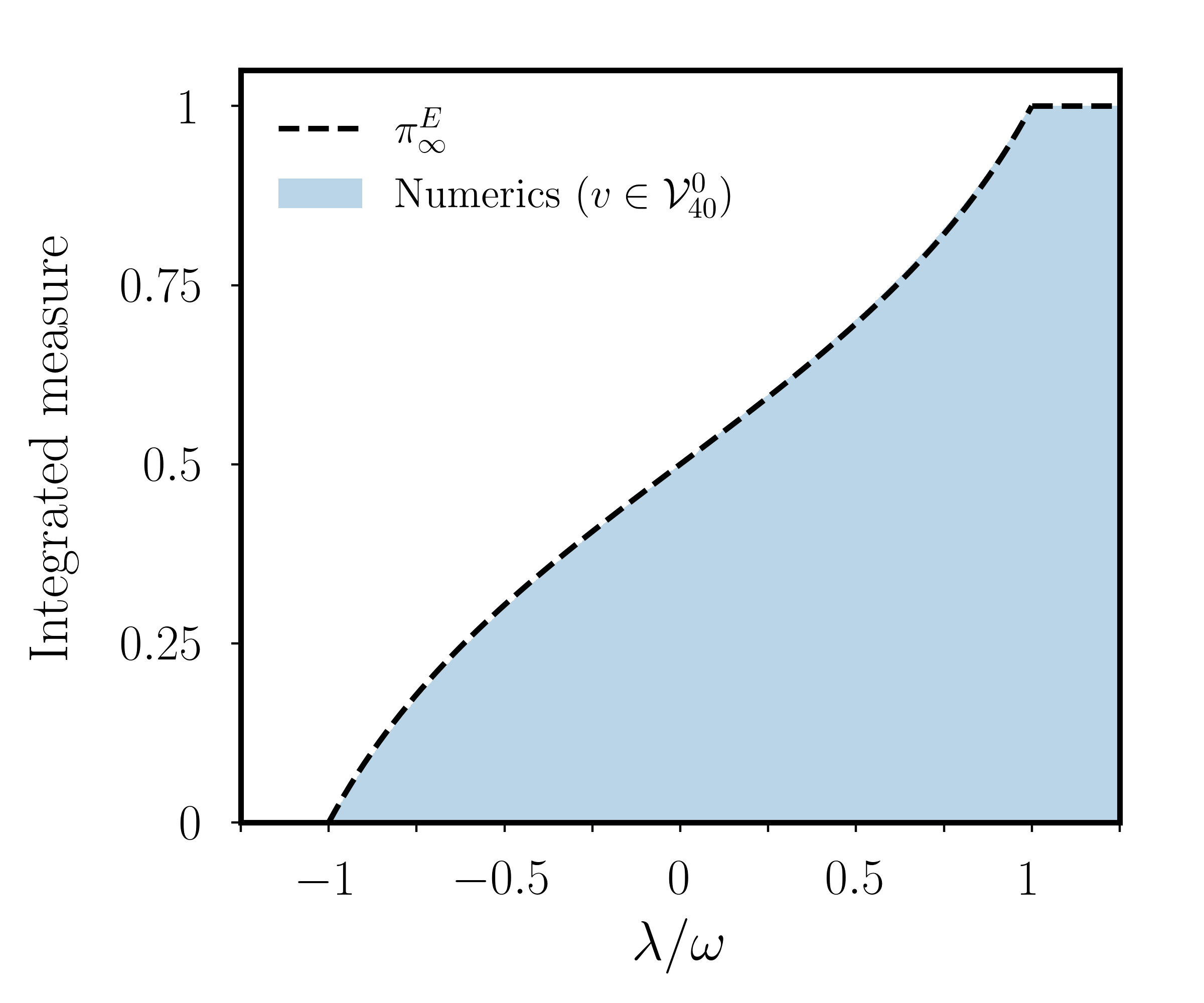} \\
        \end{tabular}
}
\qquad
\subfloat[Axisymmetric ellipsoid $A_1=A_2=1$, $A_3=0.01$ with $\theta=\pi/4$]{
    \centering
    \begin{tabular}{cc}
        \includegraphics[width=0.46\textwidth]{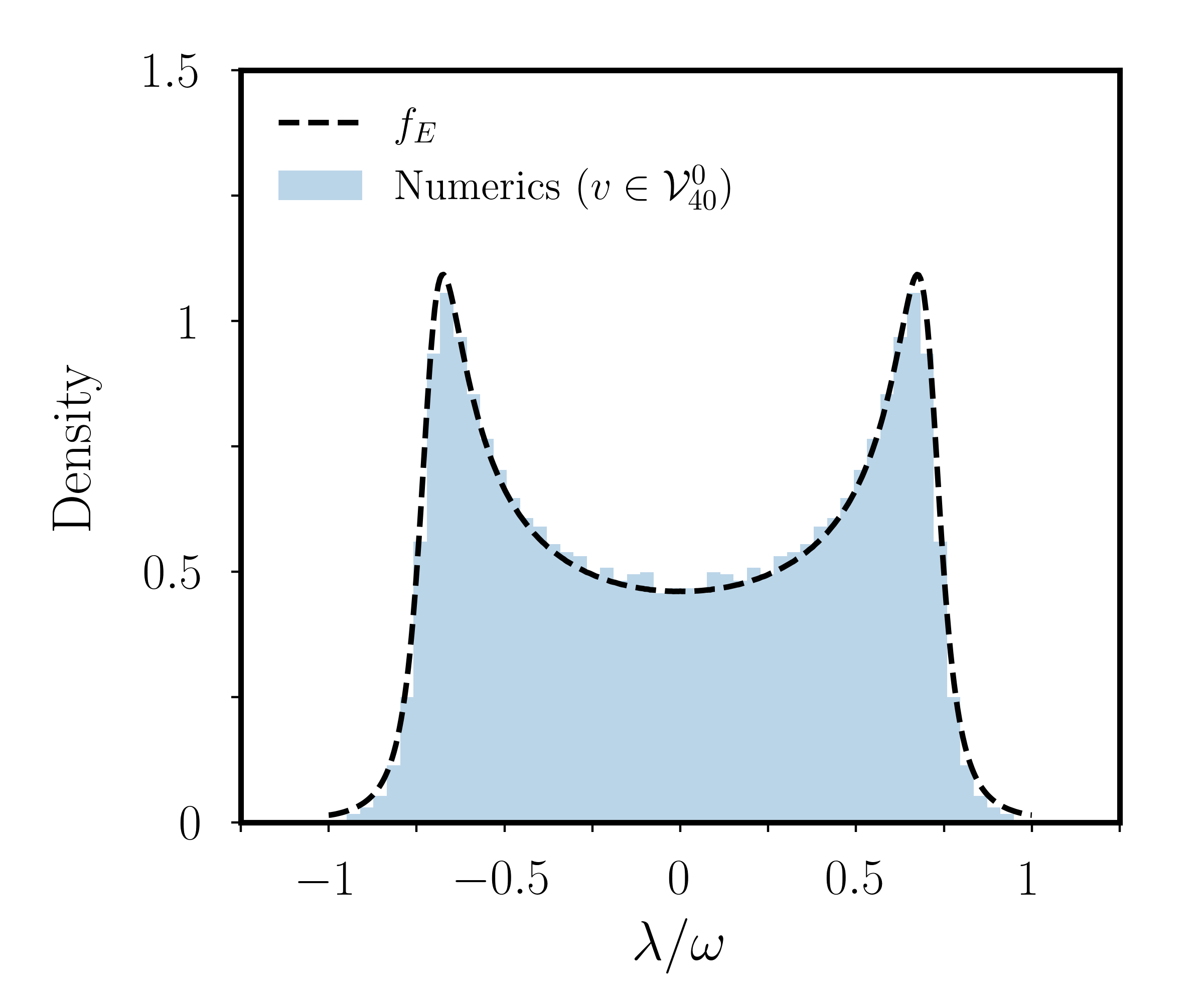}  &
        \includegraphics[width=0.46\textwidth]{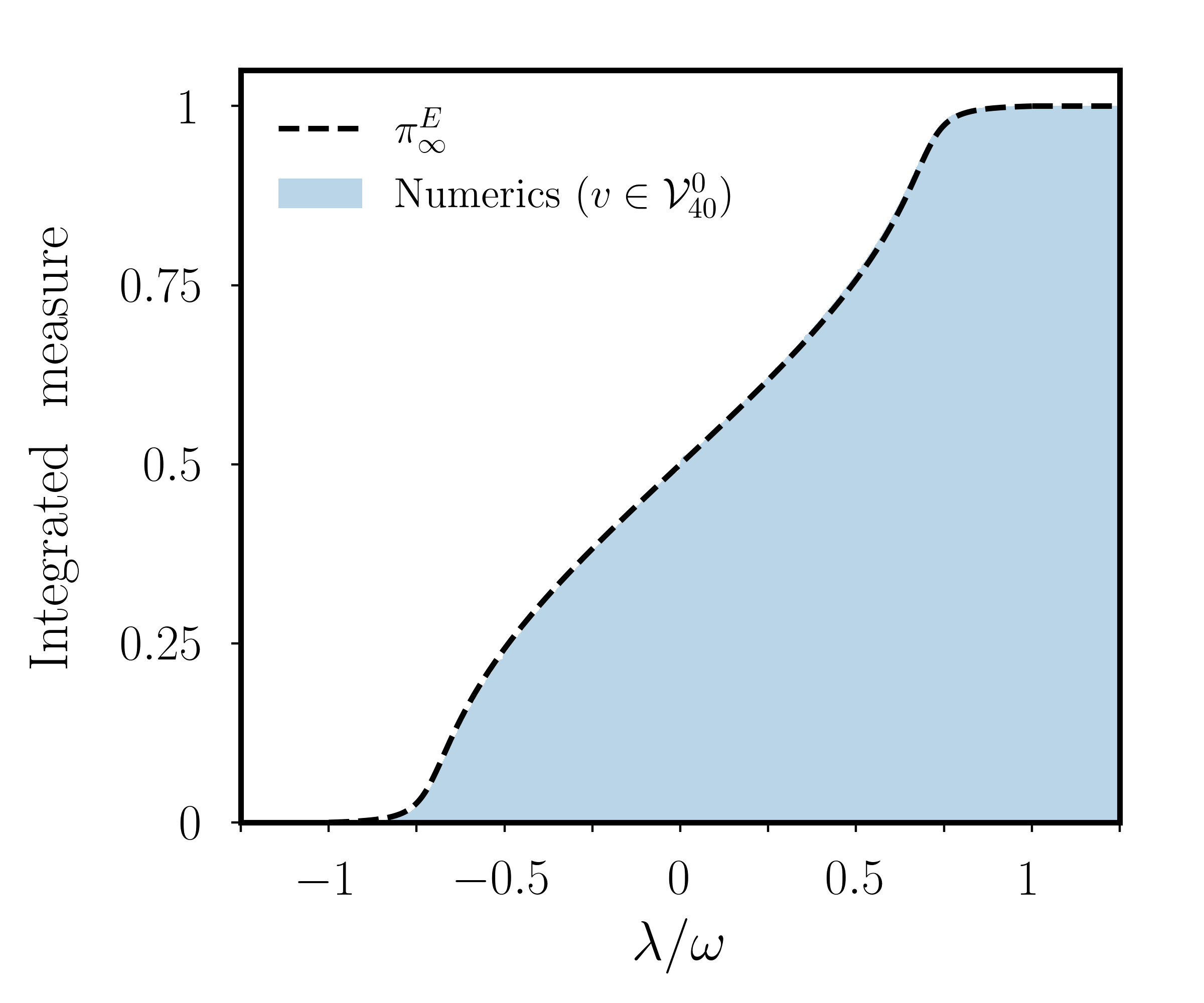} \\
    \end{tabular}}
    \caption{Comparison between numerical computations of the full inertial-mode spectrum for $v \in \mathcal{V}_{40}^0$ (i.e. considering $N=23~780$ eigenvalues) and Theorem \ref{theo:dens} for non-integrable configurations with $\vec{\Omega} = (\sin \theta, 0, \cos \theta)$.
    \emph{Left panels}: Probability densities (the $20$ degenerate eigenvalues $\lambda=0$ have been removed from the histograms). \emph{Right panels}: Integrated probability measures.}
    \label{fig:numerics}
\end{figure}

An excellent quantitative agreement is found between formula \eqref{eq:fEintegrable} and prior computations, for both the round ball (figure \ref{fig:intro}b) and other integrable cases in axisymmetric ellipsoids (not shown).
Since formula \eqref{eq:fEintegrable} is not valid when $\vec{\Omega}$ is not an axis of revolution of the geometry, it remains to compare Theorem \ref{theo:dens} with computations of the inertial-mode spectrum in such cases. 

To do so, every complex-valued eigenvector $v \in \mathcal{V}_n ^0$ of equation \eqref{equ:NSL} is sought using a polynomial expansion in ellipsoids as
\begin{subequations}
\label{eq:expandpolynomial}
\begin{equation}
    v = \sum_{j=1}^N \alpha_j \mathfrak{e}_j, \quad N = n(n+1)(2n+7)/{6},
    \tag{\theequation a--b}
\end{equation}
\end{subequations}
where $\{\mathfrak{e}_j\}_{j\leq N}$ are real-valued basis elements of $\mathcal{V}_n ^0$, and $\vec{\alpha} = (\alpha_1, \dots, \alpha_B)$ contains the complex-valued coefficients of the eigenvector in the chosen basis.
Different algorithms have been proposed to explicitly construct polynomial bases of $\mathcal{V}_n ^0$, which all involve non-orthogonal basis elements in ellipsoids \cite{lebovitz1989stability,wu2011high,ivers2017enumeration}. 
To determine the modal coefficients in expansion \eqref{eq:expandpolynomial}, we employ a projection method \cite{vidal2024igw}.
We substitute expansion \eqref{eq:expandpolynomial} in equation \eqref{equ:NSL} and, then, we project the resulting equation onto every basis element $\mathfrak{e}_i$ such that
\begin{equation*}
    i \lambda \sum_{j=1}^N \alpha_j \int_E \langle \mathfrak{e}_i \mid \mathfrak{e}_j \rangle \, dx_1 dx_2 dx_3 = -\sum_{j=1}^N \alpha_j \int_E \langle \mathfrak{e}_i \mid \vec{\Omega} \wedge \mathfrak{e}_j \rangle \, dx_1 dx_2 dx_3,
\end{equation*}
where the volume integration is performed analytically in ellipsoids \cite{lebovitz1989stability}.
This procedure gives an exact system of equations that can be written in the form of a generalized eigenvalue problem as $A \vec{\alpha} = i \lambda B \vec{\alpha}$, where $[A,B]$ are $N\times N$ real-valued matrices.
The matrix B is definite positive but, since the basis elements are not mutually orthogonal in ellipsoids, it is not diagonal.
Further details about the mathematical algorithm are given in \cite{vidal2024igw}.
For the numerics, it turns out that the two matrices $[A,B]$ usually become ill-conditioned  when $n \geq 20$ (as also found in prior numerical studies in ellipsoids \cite{gerick2020pressure,vidal2021kinematic}). 
Thus, we have numerically solved the eigenvalue problem using extended floating-point precision.

As shown in figure \ref{fig:numerics}, a very good quantitative agreement is obtained between the asymptotic measure $\pi_\infty^E$ given by Theorem \ref{theo:dens} and the numerics for arbitrary cases. 
Therefore, Theorem \ref{theo:dens} can be used to obtain the asymptotic behaviour of the inertial-mode spectrum for any configuration.

\section{Ray dynamics}
As usual the classical (ray) dynamics gives a good way to approach the eigenmodes of large degrees. 
Hence, we show below a "classical version" of Theorem \ref{theo:A}.

We define the ray dynamics associated with the Poincar\'e operator. 
We consider the Hamiltonians $\pm \Lambda_1(\xi)$ defined in Section \ref{sec:B}. 
The trajectories inside $B$ are lines because the Hamiltonians are independent of $x$. 
There is a reflection when a trajectory hits the boundary; we re-start there with the Hamiltonian field associated with the opposite eigenvalue.
The law of reflection is not trivial.
Let us consider, at energy $\gl$, the dispersion relation $\Sigma$ that is obtained as the zero set of the determinant of the symbol of $P-\gl$: we have $\Sigma = \{ \gl  \lVert \xi \rVert^2 -\go \xi_3^2 = 0 \}$.
We assume that a trajectory $(x_-(t),\xi_-)$ of the Hamiltonian $\Lambda _1$ hits the boundary at the point $m$ with the normal $\vec{n}$ to it.
Then, the reflected trajectory is $(x_+(t),\xi_+)$ where $\xi_+ $ is determined by $\xi_+\ne \xi_-$,
$\xi_+ \in \Sigma $ and $(\xi_+-\xi_-)(\vec{n})=0$ (recall that the $\xi$'s are linear forms on $\R^3$).  
The commutation of $P$ and $L$ implies the fact that the Poisson bracket of their principal symbols vanishes. 
Hence, $\gs (L)$ is a first integral of the motion.
On the other hand, the dispersion relation implies that $|\xi_3| / \lVert \xi \rVert$ is a constant along each trajectory. 
This fact can be proved directly using the law of reflection we have just described. 
This proves that the classical ray dynamics admits a constant of motion. 
If we were in dimension $2$, we would have an integrable Hamiltonian system.

\section{Perspectives}
\label{sec:persp}
The following related questions are of some interest from a mathematical viewpoint. 
\begin{enumerate}

\item An inverse spectral problem: {does the measure $\pi_\infty $ determine the pair $(E,\vec{\Omega} )$ up to dilation of $E$}?

\item {What is the asymptotic behaviour of the eigenvalues when the ellipsoid degenerates to a one-dimensional or two-dimensional limit}? 

\item {What is the asymptotic behaviour, as $n\ra \infty$, of the first eigenvalues $\gl_j^n$ for fixed $j\geq 1$ whose limits are $-\go $}?

\item {Are the non-zero eigenvalues of $P$ of finite multiplicities and, hence, all eigenfunctions polynomials?
 In case the answer is negative or unknown, are all eigenfunctions smooth}?

\item {Are there other examples of integrable cases than ellipsoids of revolution with $\vec{\Omega }$ on the axis of revolution?}
\end{enumerate}

Finally, the mathematical analysis presented in this study could also be used to better understand the properties of other linear wave motions, such as in stably stratified fluids (i.e. having a stable density profile $\rho$ such that $\langle \nabla \rho \mid \vec{g} \rangle > 0$ for incompressible fluids, where $\vec{g}$ is the gravity field).
Indeed, there is a strong analogy between uniformly rotating fluids and stably stratified fluids \cite{veronis1970analogy}, and the analogue of $\omega$ for stratified flows is called the Brunt-V\"ais\"al\"a frequency $N$.
The buoyancy force can sustain internal (gravity) waves in stratified fluids \cite{lighthill1978waves}, which usually exhibit singular spatial structures called attractors in bounded domains \cite{dintrans1999gravito,CdV-b}.
However, if we assume that $N$ and $\vec{g}$ are spatially uniform in the fluid, then the pressure associated with a divergenceless velocity field obeys a mathematical equation that is similar to Poincar\'e problem \eqref{eq:poincareeq}. 
Therefore, under these assumptions, internal (gravity) waves will admit polynomial solutions in two-dimensional elliptic domains \cite{CdV-Li} and in ellipsoids \cite{vidal2024igw,CdV-Vi-2}.

\appendix
\section{The bicharacteristic flow of $\Lambda $ and the geodesic flow on $S^3$}
\label{ss:bic}

Let us introduce the maps $j_{\pm}:B \ra S^3 \subset \R^3_x\times \R_z$
defined by
$j_\pm (x)=(x,z=\pm \sqrt{1-r^2} ) $ with $r=\lVert x \rVert$.
We have
\[h= j_\pm ^\star (dx_1^2+dx_2^2 + dx_3^2 + dz^2 )= dx_1^2+dx_2^2 + dx_3^2 + \frac{r^2}{1-r^2} dr^2. \]
Because of the invariance by the rotations in $\R^3$, we can restrict the computation of the dual metric to the points $(x_1,0,0)$ where 
\begin{equation*}
    h= \frac{1}{1-x_1^2} dx_1^2+dx_2^2 + dx_3^2.
\end{equation*}
Hence, the Hamiltonian $h^\star = (1-x_1^2)\xi_1^2+\xi_2^2+\xi_3^2$ is equal to $\Lambda (x_1,0,0;\xi)$.
This proves that the Hamiltonian flow $\phi_t$  of $\Lambda $ is the projection of the geodesic flow of $S^3$.
The geodesics of $S^3$ (which are not meridians) project onto ellipses in $\bar{B}$ tangent to $\pa B $ at two antipodal points. 
This part of the flow is simply periodic of period $2\pi $. 
The meridian great circles project onto a diameter, hitting the boundary at characteristic points of $\pa B $. 

\section{A simple Tauberian Theorem}
\begin{lemm} \label{lemm:tauber}
Let $Z(t)=\sum _{n\in \Z }a_n \mathrm{e}^{it(n+3/2)} $ be an antiperiodic distribution on $\R$. 
Assume that the singular support of $Z$ is $2\pi \Z$ and that,
for a function $\rho \in C_o^\infty (]-2\pi, 2\pi [ )$ with $\rho \equiv 1$ near $t=0$, the Fourier transform of $\rho Z $ is 
\begin{equation*}
\widehat{\rho Z }(\tau ) =b(\tau).
\end{equation*}
Then, we have $a_n = b(n+3/2) + O (n^{-\infty})$.
\end{lemm}

\begin{proof}
The asymptotics of the Fourier transform $\widehat{\rho Z }(\tau)$ is independent of $\rho $: changing $\rho$ modifies
$\rho Z $ by a smooth compactly supported function whose series has a rapid decay. 
We can choose $\rho=\rho_0 $ so that
$\sum_{n \in \Z} \rho_0 (t+2\pi n) =1$: assuming that $[-\pi, \pi ]\subset \mathrm{Support}(\rho)$, we take
\begin{equation*}
\rho _0 (t)= \rho (t)/\sum_{n\in \Z }\rho (t + 2k\pi ).
\end{equation*}
The Fourier coefficients of $Z$ are
\begin{align*}
a_n &= \frac{1}{2\pi }\int _0^{2\pi} \mathrm{e}^{-it (n+3/2)} Z (t) \, dt =
\frac{1}{2\pi }\int_\R \mathrm{e}^{-it (n+3/2)} \rho_0 (t) Z (t) \, dt, \\
{} &= b(n+3/2)+ O\left( n^{-\infty} \right). \qedhere
\end{align*}
\end{proof}

\section{The boundary \OPD~ calculus of Boutet de Monvel}
\label{app:lbm}
Boutet de Monvel developed in \cite{Bo} a pseudo-differential calculus on manifold with boundaries, and Grubb \cite{Gr} gave a more developed presentation. 
Let us assume that $X$ is a compact manifold with a smooth boundary.
The \BOPD s are operators on functions on $X$ of the form $u\ra Au:= P\tilde{u}+G\tilde{u}$  where $\tilde{u}$ is the extension of $u$ by $0$ outside
$\bar{X}$ and
\begin{enumerate}
\item $P$ is a (classical) \OPD~in some neighbourhood of $\bar{X}$ satisfying the \textit{transmission property}; this property is satisfied for differential operators and
their parametrices.
\item $G$ is an integral operator that "lives" near the boundary $\pa X$.
\end{enumerate}

More precisely, $G$ is a sum of operators involving traces on the boundary and an integral operator involving a Fourier transform at the boundary (see definition 10.8 in \cite{Gr}). 
What will be important for us is that such operators form an algebra: the \OPD ~ parts $P$  compose with the usual rules; if $\Delta $ is an elliptic operator with elliptic boundary conditions, then the parametrix of $\Delta $ is of the previous form.

The main result that we used (section \ref{ss:generalcase}) is the following \cite[proposition 10.11]{Gr}.
\begin{prop} If $\chi $ is a smooth function vanishing near the boundary and $G$ is a Green operator,
$\chi G$ and $G \chi $ are smoothing operators.
\end{prop}

Our basic example is the Leray projector $\Pi =dd^\star \Delta^{-1}$ where $\Delta=dd^\star + d^\star d$  is the Hodge-de Rham laplacian with relative boundary conditions:
if $j:\pa B \ra \bar{B}$ is the embedding, we ask that the $2$-forms $\omega $ satisfy   $j^\star \omega =0,~j^\star (d^\star \omega )=0$.
Recall that $\Delta $ is invertible: the kernel of $\Delta $ is isomorphic to the space $H^2_\mathrm{rel}(B,\pa B;\R)$ that vanishes, because $B$ is simply connected and the Poincar\'e duality. 
The Poincar\'e operator, which is a composition of Leray projectors and a wedge product by $\vec{\Omega}$, thus belongs also to this algebra. 

\section{A short review about the calculus of wave-front sets}
\label{app:wf}
For more details on this section, we can look \cite[section 2.5]{Ho-1} or \cite[section 1.3]{Du}.
Let us simply recall that, if $u$ is a Schwartz distribution on a smooth manifold $X$, one can define the wave-front set $\mathrm{WF}(u)$ of $u$ as a closed conical subset of $T^\star X\setminus 0$ whose projection onto $X$ is the singular support of $X$.

If $A:C_o^\infty (Y) \ra \mathcal{D}'(X)$ is a linear operator, we can look at the Schwartz kernel $[A]$ of $A$, which is a distribution on $X\times Y$. 
It is then natural to define 
\begin{equation*}
\mathrm{WF}'(A):= \{ (x,\xi;y, \eta )\in T^\star (X\times Y) \mid (x,\xi;y,-\eta ) \in \mathrm{WF}([A]) \}.
\end{equation*}
If $\mathrm{WF}'(A)\subset T^\star (X)\setminus 0 \times  T^\star (Y)\setminus 0$, then we have $\mathrm{WF}(Au)\subset \mathrm{WF}'(A)\circ \mathrm{WF}(u) $ for any $u\in \mathcal{D}'(Y)$.

Hörmander \cite{Ho-1} introduced what is called Fourier Integral Distributions associated with a conic Lagrangian manifold $Z$. 
Such distributions are defined as oscillatory integrals of the form
\begin{equation*}
u(x)=\int_{\R^N} \mathrm{e}^{i\phi(x,\theta )} a(x,\theta) \, |d\theta|,
\end{equation*}
where $\phi $ is "generating function" of $Z$ and $a$ is a symbol.
There exists $m\in \R$ so that, for any multi-indices $(\ga, \gb)$, there exists $C_{\ga,\gb}>0$ so that
\begin{equation*}
|D_x^\ga D_\theta^\gb |(x,\theta) \leq C_{\ga,\gb}|\theta |^{m-|\beta |}.
\end{equation*}
If $u$ is such a distribution, then we have $\mathrm{WF}(u)\subset Z$. 


\begin{ack}
This paper is dedicated to the memory of Steve. Many thanks to Charles Epstein, Suresh Eswarathasan, Gerd Grubb, Cyril Letrouit, Johannes Sjöstrand, Andras Vasy, Bernard Valette, Jared Wunsch, Yuan Xu for answers to our questions and comments about preliminary versions. 
\end{ack}

\begin{funding}
JV received funding from European Research Council under the European Union's Horizon 2020 research and innovation programme (\textsc{theia} project, grant agreement no. 847433). 
The numerical computations presented in this paper were performed using the \textsc{gricad} infrastructure (\url{https://gricad.univ-grenoble-alpes.fr}), which is supported by Grenoble research communities.
\end{funding}

\bibliography{biblio}
\bibliographystyle{emss}







\end{document}